\documentclass[11pt]{amsart}
\usepackage{amsmath,amssymb,amsthm}
\usepackage[margin=1in]{geometry}

\newtheorem{theorem}{Theorem}[section]
\newtheorem*{theorem*}{Theorem}
\newtheorem{lemma}[theorem]{Lemma}
\newtheorem*{lemma*}{Lemma}
\newtheorem{corollary}[theorem]{Corollary}
\newtheorem*{corollary*}{Corollary}
\newtheorem{proposition}[theorem]{Proposition}
\newtheorem*{proposition*}{Proposition}
\theoremstyle{definition}
\newtheorem{definition}[theorem]{Definition}
\newtheorem*{definition*}{Definition}

\newtheorem*{example*}{Example}
\newtheorem{remark}[theorem]{Remark}
\newtheorem*{remark*}{Remark}

\newtheorem*{notation*}{Notation}

\begin{document}

\title{Irreducibility of quadratic hulls of osculating varieties of closed orbits}
\author{César Massri}
\begin{abstract}
Let $G$ be a connected complex semisimple algebraic group, let $V_\lambda$ be
an irreducible $G$-module of highest weight $\lambda$, and let
\[
    X_{d\lambda}=G\cdot[v_\lambda^d]
    \subset \mathbb P V_{d\lambda}
\]
be the closed orbit of a highest weight vector. We study the irreducibility of
the zero locus of the quadratic equations containing the $p$-osculating variety
$T^pX_{d\lambda}$. We prove that, if $d>4p$, its quadratic hull coincides with
the osculating variety itself:
\[
    \mathcal Q(T^pX_{d\lambda})=T^pX_{d\lambda}.
\]
In particular, the quadratic hull is irreducible. Thus, for every fixed $p$,
sufficiently high Cartan powers of a closed highest weight orbit have
$p$-osculating varieties determined set-theoretically by their quadratic
equations.
\end{abstract}


\maketitle

\section*{Introduction}

Let $Y\subset\mathbb P W$ be a projective variety. We denote by
\[
    \mathcal Q(Y)=Z(I(Y)_2)
\]
its quadratic hull, that is, the zero locus of all quadratic equations
containing $Y$. Throughout, zero loci are considered set-theoretically, or
equivalently with their reduced structure. In general, $\mathcal Q(Y)$ may be
strictly larger than $Y$ and need not be irreducible. In this article we study this problem for
osculating varieties of closed orbits in irreducible representations of
complex semisimple groups.

The motivation comes from the case of binary forms. Let $V=\mathbb C^2$ and
recall that
\[
    S^2(S^rV)
    =
    \bigoplus_{m=0}^{\lfloor r/2\rfloor} S^{2r-4m}V.
\]
In \cite{Massri2013}, for every $m$ we considered the variety denoted there by
$M_m$, defined by the quadratic equations corresponding to the irreducible
summand $S^{2r-4m}V$. To avoid confusion with the notation used below, we
denote this variety here by $Z_m$. More precisely, if
\[
    \pi_m:S^2(S^rV)\longrightarrow S^{2r-4m}V
\]
is the projection, then
\[
    Z_m
    =
    \{x\in\mathbb P S^rV\mid \pi_m(x\cdot x)=0\}.
\]
These varieties are not irreducible in general. For instance, the case
$r=8$, $m=3$ was studied by B\"ohning \cite{Bohning}; the corresponding
variety $Z_3\subset\mathbb P S^8V$ has two irreducible components.

Nevertheless, the computations in \cite{Massri2013} suggest a different
phenomenon when several irreducible summands of $S^2(S^rV)$ are considered
simultaneously. Let
\[
    c_r\subset\mathbb P S^rV
\]
be the Veronese curve and let $T^pc_r$ denote its $p$-osculating variety.
Recall that
\[
    I(T^pc_r)_2
    =
    \bigoplus_{m=p+1}^{\lfloor r/2\rfloor}
    \bigl(S^{2r-4m}V\bigr)^*.
\]
It follows that
\[
    \mathcal Q(T^pc_r)
    =
    Z_{p+1}\cap\cdots\cap Z_{\lfloor r/2\rfloor}.
\]
Although the individual varieties $Z_m$ may have extra irreducible
components, the computations in \cite{Massri2013} suggested that, in a
suitable range, this intersection coincides with the osculating variety
$T^pc_r$ itself.

Outside the range where this equality holds, the behavior may be different.
For example, consider the third osculating variety of the rational normal curve
$c_{10}\subset\mathbb P^{10}$. In this case
\[
    I(T^3c_{10})_2
    =
    (S^4V)^*\oplus\mathbb C,
\]
so $\mathcal Q(T^3c_{10})$ is defined by six quadrics. The computation in
\cite{Massri2013} gives
\[
    \dim \mathcal Q(T^3c_{10})=4.
\]
Hence these six quadrics have codimension $6$ in $\mathbb P^{10}$ and form a
regular sequence. Thus $\mathcal Q(T^3c_{10})$ is a complete intersection of
six quadrics; in particular, it is equidimensional of dimension $4$ and
\[
    \deg \mathcal Q(T^3c_{10})=2^6=64.
\]
On the other hand,
\[
    \dim T^3c_{10}=4,
    \qquad
    \deg T^3c_{10}=28.
\]
If $T^3c_{10}$ were the only irreducible component of maximal dimension of
$\mathcal Q(T^3c_{10})$, then
\[
    \deg \mathcal Q(T^3c_{10})
    =
    k\,\deg T^3c_{10}
\]
for some positive integer $k$. Since $28$ does not divide $64$, the quadratic
hull has another irreducible component of dimension $4$. Thus even in the
binary case the quadratic hull of an osculating variety may be reducible.

These examples show that neither $G$-invariance nor the geometric origin of
the quadratic equations is enough, by itself, to guarantee irreducibility.
The point of the result below is that the particular family of quadrics
determined by the PBW filtration has a stable behavior when the highest weight
is sufficiently scaled.

We now turn to the general setting. Let $G$ be a connected complex semisimple
algebraic group with Lie algebra $\mathfrak g$, let $V_\lambda$ be an
irreducible $G$-module of highest weight $\lambda$, and let
$v_\lambda\in V_\lambda$ be a highest weight vector. For every $d\geq 1$, the
representation $V_{d\lambda}$ is the Cartan component of $S^d(V_\lambda)$ and
its highest weight vector may be identified with $v_\lambda^d$. We denote by
\[
    X_{d\lambda}
    =
    G\cdot[v_\lambda^d]
    \subset
    \mathbb P V_{d\lambda}
\]
the corresponding closed orbit. Thus the underlying homogeneous variety
$G/P_\lambda$ is fixed, while its projective embedding is replaced by the
$d$-th Cartan power.

The PBW filtration gives a natural description of the osculating spaces of
$X_{d\lambda}$. More precisely, let
\[
    E_p(d\lambda)
    =
    U_{\leq p}(\mathfrak g)v_\lambda^d
    \subset
    V_{d\lambda},
\]
where $U_{\leq p}(\mathfrak g)$ is the $p$-th term of the PBW filtration of
the universal enveloping algebra. Then $\mathbb P E_p(d\lambda)$ is the
$p$-osculating space of $X_{d\lambda}$ at $[v_\lambda^d]$, and
\[
    T^pX_{d\lambda}
    =
    \overline{G\cdot\mathbb P(E_p(d\lambda))}.
\]
We study the quadratic hull $\mathcal Q(T^pX_{d\lambda})$. There is always an
inclusion
\[
    T^pX_{d\lambda}
    \subseteq
    \mathcal Q(T^pX_{d\lambda}),
\]
but in general equality need not hold. Our main result shows that equality
holds in a stable range.

\medskip

\noindent
\textbf{Theorem.}
\textit{
Let $G$ be a connected complex semisimple algebraic group, let $V_\lambda$ be
an irreducible $G$-module of highest weight $\lambda$, and let $p\geq0$ and
$d\geq1$ be integers. If
\[
    d>4p,
\]
then
\[
    \mathcal Q(T^pX_{d\lambda})
    =
    T^pX_{d\lambda}.
\]
In particular, the quadratic hull $\mathcal Q(T^pX_{d\lambda})$ is irreducible.
}

\medskip

The theorem gives a positive answer to the irreducibility problem for the
family of quadratic equations naturally determined by the PBW filtration. In
fact, the result is stronger: it identifies their zero locus with the
corresponding osculating variety.

For $\mathfrak g=\mathfrak{sl}_2(\mathbb C)$, this recovers the phenomenon
suggested by the computations in \cite{Massri2013}. The main technical
ingredient in this case is a result of Weyman on binary forms
\cite{WeymanGordan}. The general proof uses the binary case to obtain a
factorization result and then relates this factorization with the geometry of
the closed orbit and the PBW filtration.

This article is divided into four sections. In Section~1 we recall the relation
between the PBW filtration and the osculating spaces of a closed orbit, and we
describe the quadratic equations containing these varieties. In Section~2 we
study the case of binary forms and recall the result of Weyman used in the
proof. Section~3 contains the proof of the main theorem. Finally, in Section~4
we study the behavior outside the stable range: we give a binary example whose
quadratic hull is irreducible but strictly larger than the osculating variety,
discuss a reducible example related to restricted chordal varieties of
Grassmannians, and finish with a related example in which the cubic equations
of the trifocal variety, after a second Veronese embedding, give rise to a
reducible quadratic hull.

\section{PBW filtrations and osculating varieties}

Let $G$ be a connected complex semisimple algebraic group with Lie algebra
$\mathfrak g$, and let $V_\nu$ be the irreducible $G$-module of highest
weight $\nu$. Fix a highest weight vector $v_\nu\in V_\nu$, and consider
the closed orbit
\[
X_\nu=G\cdot[v_\nu]\subset \mathbb P(V_\nu).
\]

Recall that the PBW filtration on the universal enveloping algebra
$U(\mathfrak g)$ is given by
\[
U_{\leq p}(\mathfrak g)
=
\operatorname{span}
\left\{
\xi_1\cdots\xi_k:
\xi_i\in\mathfrak g,\ 0\leq k\leq p
\right\},
\]
where $U_{\leq 0}(\mathfrak g)=\mathbb C$. Thus
\[
U_{\leq 0}(\mathfrak g)
\subset
U_{\leq 1}(\mathfrak g)
\subset
U_{\leq 2}(\mathfrak g)
\subset\cdots
\subset U(\mathfrak g).
\]
For $p\geq0$, set
\[
E_p(\nu):=U_{\leq p}(\mathfrak g)v_\nu\subset V_\nu.
\]

We denote by $T_x^pX_\nu$ the $p$-th osculating
space of $X_\nu$ at $x$. Since $X_\nu$ is homogeneous, we write
\[
T^pX_\nu
=
\overline{G\cdot T^p_{[v_\nu]}X_\nu}
\]
for the corresponding osculating variety.

The relation between the PBW filtration and osculating spaces is
well known; see \cite[Proposition~2.3]{LandsbergManivel}. We recall
it here for completeness.

\begin{proposition}\label{prop:pbw-osculating}
The $p$-th osculating space of $X_\nu$ at the highest weight point is
\[
T^p_{[v_\nu]}X_\nu
=
\mathbb P(E_p(\nu)).
\]
Consequently,
\[
T^pX_\nu
=
\overline{G\cdot\mathbb P(E_p(\nu))}.
\]
\end{proposition}

\begin{proof}
Let $P_\nu$ be the stabilizer of $[v_\nu]$, with Lie algebra
$\mathfrak p_\nu$, and let $\mathfrak n_\nu^-$ be the opposite
nilradical. Then
\[
\mathfrak g=\mathfrak n_\nu^-\oplus\mathfrak p_\nu.
\]
Choose an ordered basis of $\mathfrak g$ consisting first of a basis
of $\mathfrak n_\nu^-$ and then of a basis of $\mathfrak p_\nu$.
By the Poincar\'e--Birkhoff--Witt theorem
\cite[\S17.3]{Humphreys},
$U_{\leq p}(\mathfrak g)$ is spanned by elements
$u_-u_+$, where
\[
u_-\in U_{\leq a}(\mathfrak n_\nu^-),\qquad
u_+\in U_{\leq b}(\mathfrak p_\nu),\qquad a+b\leq p.
\]
Since $\mathfrak p_\nu$ preserves the line $\mathbb Cv_\nu$, we have
\[
U(\mathfrak p_\nu)v_\nu\subset\mathbb Cv_\nu.
\]
Hence
\[
U_{\leq p}(\mathfrak g)v_\nu
=
U_{\leq p}(\mathfrak n_\nu^-)v_\nu.
\]

Choose an ordered basis $\xi_1,\ldots,\xi_r$ of $\mathfrak n_\nu^-$.
Since the natural map
$\mathfrak n_\nu^-\to\mathfrak g/\mathfrak p_\nu$ is an isomorphism, the map
\[
(t_1,\ldots,t_r)\longmapsto
\left[
\exp(t_1\xi_1)\cdots\exp(t_r\xi_r)v_\nu
\right]
\]
has invertible differential at the origin and gives local coordinates on
$X_\nu$ around $[v_\nu]$.
Using
\[
\widetilde\phi(t_1,\ldots,t_r)
=
\exp(t_1\xi_1)\cdots\exp(t_r\xi_r)v_\nu
\]
as a local lift, we have
\[
\frac{\partial^{|\alpha|}\widetilde\phi}
{\partial t_1^{\alpha_1}\cdots\partial t_r^{\alpha_r}}(0)
=
\xi_1^{\alpha_1}\cdots\xi_r^{\alpha_r}v_\nu.
\]
By the PBW theorem, the ordered monomials
\[
\xi_1^{a_1}\cdots\xi_r^{a_r},
\qquad
a_1+\cdots+a_r\leq p,
\]
span $U_{\leq p}(\mathfrak n_\nu^-)$. Hence the derivatives of
order at most $p$ span
\[
U_{\leq p}(\mathfrak n_\nu^-)v_\nu
=
E_p(\nu).
\]
It follows that
\[
T^p_{[v_\nu]}X_\nu
=
\mathbb P(E_p(\nu)).
\]

Finally, by $G$-equivariance,
\[
T^p_{g[v_\nu]}X_\nu
=
gT^p_{[v_\nu]}X_\nu,
\]
and therefore
\[
T^pX_\nu
=
\overline{G\cdot\mathbb P(E_p(\nu))}.
\]
\end{proof}

We now describe the quadratic equations of this variety.
For a subspace $W\subset S^2(V_\nu)$, we denote by
\[
W^\circ
=
\{q\in S^2(V_\nu^*):q(w)=0\text{ for every }w\in W\}
\]
its annihilator. Set
\[
F_p(\nu):=
U(\mathfrak g)\cdot S^2(E_p(\nu))
\subset S^2(V_\nu).
\]

\begin{proposition}\label{prop:quadrics}
One has
\[
I(T^pX_\nu)_2
=
F_p(\nu)^\circ
\subset S^2(V_\nu^*).
\]
\end{proposition}

\begin{proof}
A quadratic form $q\in S^2(V_\nu^*)$ vanishes on $T^pX_\nu$ if and
only if it vanishes on the affine cone $gE_p(\nu)$ for every $g\in G$.
By polarization, a quadratic form vanishes identically on a linear
subspace $W\subset V_\nu$ if and only if it annihilates
$S^2(W)\subset S^2(V_\nu)$. Therefore, $q$ vanishes on $T^pX_\nu$
if and only if it annihilates
\[
g\cdot S^2(E_p(\nu))
\]
for every $g\in G$.

Since $G$ is connected,
\[
\operatorname{span}
\{g\cdot S^2(E_p(\nu)):g\in G\}
=
U(\mathfrak g)\cdot S^2(E_p(\nu))
=
F_p(\nu).
\]
Hence $q$ annihilates $F_p(\nu)$.
Thus
\[
I(T^pX_\nu)_2=F_p(\nu)^\circ.
\]
\end{proof}

\begin{definition}
We set
\[
M_p(G,V_\nu)
:=
\mathcal Q(T^pX_\nu)
=
Z\bigl(I(T^pX_\nu)_2\bigr)
\subset \mathbb P(V_\nu).
\]
\end{definition}

Here, for $x\in V_\nu$, we write $x\cdot x$ for its symmetric product
in $S^2(V_\nu)$. By Proposition~\ref{prop:quadrics} and finite-dimensional
duality,
\[
[x]\in M_p(G,V_\nu)
\quad\Longleftrightarrow\quad
q(x)=0\text{ for every }q\in F_p(\nu)^\circ
\quad\Longleftrightarrow\quad
x\cdot x\in(F_p(\nu)^\circ)^\circ=F_p(\nu).
\]
Thus
\[
M_p(G,V_\nu)
=
\left\{
[x]\in\mathbb P(V_\nu)
\;:\;
x\cdot x\in F_p(\nu)
\right\}.
\]

By construction,
\[
T^pX_\nu
\subseteq
M_p(G,V_\nu).
\]
Thus the problem considered in this article is to determine when the
reverse inclusion holds.

\section{The binary case}

Let $V=\mathbb C^2$ and $G=SL(V)$. For $r\geq 1$, consider the
irreducible representation $S^rV$ and its closed orbit
\[
c_r
=
\{[\ell^r]:0\neq\ell\in V\}
\subset
\mathbb P(S^rV),
\]
the rational normal curve.

We write $E_p(r)$ and $F_p(r)$ for the spaces $E_p(\nu)$ and
$F_p(\nu)$ introduced in Section~1 in the case $V_\nu=S^rV$.

Choose a basis $u,v$ of $V$ and let
\[
H=
\begin{pmatrix}
1&0\\
0&-1
\end{pmatrix},
\qquad
Y=
\begin{pmatrix}
0&0\\
1&0
\end{pmatrix}
\in\mathfrak{sl}_2.
\]
Thus $u^r$ is a highest weight vector and $Yu=v$. Set
\[
x_i=\frac{Y^iu^r}{i!},
\qquad
0\leq i\leq r.
\]
Then $x_i$ has weight $r-2i$, and $x_i$ is a non-zero scalar
multiple of $u^{r-i}v^i$.

We first recall the geometric description of the osculating
varieties of $c_r$.

\begin{proposition}\label{prop:binary-osculating}
For $0\leq p\leq r$,
\[
T^pc_r
=
\left\{
[\ell^{r-p}h]:
0\neq\ell\in V,\ 0\neq h\in S^pV
\right\}.
\]
In particular, $T^pc_r$ is the variety of binary forms of degree $r$
having a root of multiplicity at least $r-p$.
\end{proposition}

\begin{proof}
By Proposition~\ref{prop:pbw-osculating},
\[
T^p_{[u^r]}c_r
=
\mathbb P(E_p(r)).
\]
The opposite nilradical of $\mathfrak{sl}_2$ is generated by $Y$.
Hence
\[
E_p(r)
=
\operatorname{span}\{x_0,\ldots,x_p\}
=
u^{r-p}S^pV.
\]
Since $SL(V)$ acts transitively on $\mathbb P(V)$, it follows that
\[
T^pc_r
=
\left\{
[\ell^{r-p}h]:
0\neq\ell\in V,\ 0\neq h\in S^pV
\right\}.
\]
Indeed, the set on the right is closed, since it is the image of the
projective morphism
\[
\mathbb P(V)\times\mathbb P(S^pV)
\longrightarrow
\mathbb P(S^rV),
\qquad
([\ell],[h])
\longmapsto
[\ell^{r-p}h].
\]
\end{proof}

We now determine $F_p(r)$. Recall the Clebsch--Gordan decomposition
\[
S^2(S^rV)
=
\bigoplus_{m=0}^{\lfloor r/2\rfloor}
S^{2r-4m}V;
\]
see, for instance, \cite{Massri2013}. In particular, this
decomposition is multiplicity free. For every
$0\leq m\leq\lfloor r/2\rfloor$, let
\[
\pi_m:
S^2(S^rV)
\longrightarrow
S^{2r-4m}V
\]
be a non-zero $SL(V)$-equivariant projection.

Fix a highest weight vector $w_m$ of $S^{2r-4m}V$. For
$z\in S^2(S^rV)$, let $q_m(z)$ be the coefficient of $w_m$ in the
weight-$(2r-4m)$ component of $\pi_m(z)$.

The following elementary computation will be useful. It is the
coefficient computation carried out in \cite[Corollary~2.2]{Massri2013},
rewritten in the present notation.

\begin{lemma}\label{lem:binary-coefficients}
For every $0\leq m\leq\lfloor r/2\rfloor$, there exists a non-zero
constant $\lambda_m$ such that, for $0\leq i,j\leq r$,
\[
q_m(x_i\cdot x_j)
=
\begin{cases}
(-1)^i\binom{2m}{i}\lambda_m,
    & \text{if }i+j=2m,\\[4pt]
0,  & \text{otherwise}.
\end{cases}
\]
\end{lemma}

\begin{proof}
The vector $x_i\cdot x_j$ has weight
\[
2r-2(i+j),
\]
whereas $w_m$ has weight $2r-4m$. Since $\pi_m$ preserves weights,
$q_m(x_i\cdot x_j)$ can be non-zero only if
\[
i+j=2m.
\]

Set
\[
a_i=q_m(x_i\cdot x_{2m-i}).
\]
For $0\leq i<2m$, the vector
$x_i\cdot x_{2m-i-1}$ has weight $2r-4m+2$, which is larger than
the highest weight of $S^{2r-4m}V$. Hence
\[
\pi_m(x_i\cdot x_{2m-i-1})=0.
\]
Using the equivariance of $\pi_m$ and
\[
Yx_i=(i+1)x_{i+1},
\]
we obtain
\[
0
=
\pi_m\bigl(Y(x_i\cdot x_{2m-i-1})\bigr).
\]
Taking the coefficient of $w_m$ gives
\[
(i+1)a_{i+1}+(2m-i)a_i=0.
\]
Therefore
\[
a_i=(-1)^i\binom{2m}{i}a_0.
\]

Since $\pi_m$ is non-zero and $S^{2r-4m}V$ is irreducible,
$\pi_m$ is surjective. In particular, some vector of weight $2r-4m$
in $S^2(S^rV)$ maps to a non-zero multiple of $w_m$. Hence $q_m$ is
non-zero on the weight-$2r-4m$ subspace of $S^2(S^rV)$. This weight
space has basis
\[
x_i\cdot x_{2m-i},
\qquad
0\leq i\leq m.
\]
By the preceding recursion, all the values of $q_m$ on these basis
vectors are scalar multiples of $a_0$. It follows that $a_0\neq0$.
Taking $\lambda_m=a_0$ proves the result.
\end{proof}

We can now describe $F_p(r)$.

\begin{proposition}\label{prop:binary-F}
For $0\leq p\leq r$,
\[
F_p(r)
=
\bigoplus_{m=0}^{\min\{p,\lfloor r/2\rfloor\}}
S^{2r-4m}V.
\]
\end{proposition}

\begin{proof}
Since
\[
E_p(r)
=
\operatorname{span}\{x_0,\ldots,x_p\},
\]
the space $S^2(E_p(r))$ is spanned by
\[
x_i\cdot x_j,
\qquad
0\leq i,j\leq p.
\]

Suppose first that $m>p$. Every such vector has weight
\[
2r-2(i+j)
\geq
2r-4p
>
2r-4m.
\]
Since $2r-4m$ is the highest weight of $S^{2r-4m}V$, we have
\[
\pi_m(x_i\cdot x_j)=0
\]
for every $0\leq i,j\leq p$. By equivariance, $\pi_m$ also vanishes
on
\[
U(\mathfrak{sl}_2)\cdot S^2(E_p(r))
=
F_p(r).
\]
Thus the component $S^{2r-4m}V$ does not occur in $F_p(r)$.

Conversely, let
\[
0\leq m\leq\min\{p,\lfloor r/2\rfloor\}.
\]
Then $x_m\in E_p(r)$ and, by Lemma~\ref{lem:binary-coefficients},
\[
q_m(x_m\cdot x_m)
=
(-1)^m\binom{2m}{m}\lambda_m
\neq0.
\]
Hence
\[
\pi_m(S^2(E_p(r)))\neq0,
\]
and therefore $\pi_m(F_p(r))\neq0$.

Finally, $F_p(r)$ is an $SL(V)$-submodule of $S^2(S^rV)$.
Finite-dimensional $SL(V)$-modules are completely reducible, and the
Clebsch--Gordan decomposition above is multiplicity free. Hence
$F_p(r)$ is the direct sum of a subset of its irreducible summands.
The preceding two paragraphs show that this subset consists precisely
of the summands with
\[
0\leq m\leq\min\{p,\lfloor r/2\rfloor\}.
\]
\end{proof}

Combining this proposition with Proposition~\ref{prop:quadrics} gives the quadratic
equations of the osculating variety.

\begin{corollary}\label{cor:binary-quadrics}
For $0\leq p\leq r$,
\[
I(T^pc_r)_2
=
\bigoplus_{m=p+1}^{\lfloor r/2\rfloor}
\bigl(S^{2r-4m}V\bigr)^*.
\]
Consequently,
\[
M_p(SL(V),S^rV)
=
\left\{
[f]\in\mathbb P(S^rV):
\pi_m(f\cdot f)=0
\text{ for every }
p+1\leq m\leq\lfloor r/2\rfloor
\right\}.
\]
\end{corollary}

\begin{proof}
By Proposition~\ref{prop:quadrics},
\[
I(T^pc_r)_2
=
F_p(r)^\circ.
\]
The first assertion follows from the preceding proposition and the
Clebsch--Gordan decomposition. The second follows from the definition
\[
M_p(SL(V),S^rV)
=
\mathcal Q(T^pc_r)
\]
and from the fact that
\[
f\cdot f\in F_p(r)
\]
if and only if all its components in $S^{2r-4m}V$, $m>p$, vanish.
\end{proof}

We now recall the result on binary forms that will be used below.
Recall that the classical $2m$-th transvectant is symmetric in
its two arguments and defines a non-zero $SL(V)$-equivariant map
\[
(\,\cdot\,,\,\cdot\,)_{2m}:
S^2(S^rV)
\longrightarrow
S^{2r-4m}V;
\]
see \cite[Section~1(D)]{WeymanGordan}. Since $S^{2r-4m}V$ occurs with
multiplicity one in $S^2(S^rV)$, this map is a non-zero scalar multiple
of $\pi_m$. Consequently,
\[
(f,f)_{2m}=0
\quad\Longleftrightarrow\quad
\pi_m(f\cdot f)=0.
\]

The following consequence of Weyman's Theorem~2.15 is the
key input in the binary case.

\begin{proposition}\label{prop:weyman-binary}
Let $0\neq f\in S^rV$ and let $k\geq1$. Assume that
\[
r>4k-4.
\]
If
\[
(f,f)_{2m}=0
\qquad
\text{for every }
k\leq m\leq\lfloor r/2\rfloor,
\]
then $f$ has a root of multiplicity at least $r-k+1$.
\end{proposition}

\begin{proof}
The hypothesis $r>4k-4$ is equivalent to
\[
2k-2<\frac r2.
\]
By Proposition~\ref{prop:binary-osculating}, the locus of forms in
$\mathbb P(S^rV)$ having a root of multiplicity at least $r-k+1$ is
precisely
\[
T^{k-1}c_r.
\]
Under the inequality above, \cite[Theorem~2.15]{WeymanGordan} says
that the homogeneous ideal of this locus is generated by its
degree-two part. By Corollary~\ref{cor:binary-quadrics}, applied with
$p=k-1$, this degree-two part is
\[
I(T^{k-1}c_r)_2
=
\bigoplus_{m=k}^{\lfloor r/2\rfloor}
\bigl(S^{2r-4m}V\bigr)^*.
\]
The hypothesis
\[
(f,f)_{2m}=0,
\qquad
k\leq m\leq\lfloor r/2\rfloor,
\]
is equivalent to
\[
\pi_m(f\cdot f)=0,
\qquad
k\leq m\leq\lfloor r/2\rfloor,
\]
and hence to the vanishing at $f$ of every quadratic equation of
$T^{k-1}c_r$. Since the homogeneous ideal of $T^{k-1}c_r$ is
generated by these quadrics, it follows that
\[
[f]\in T^{k-1}c_r.
\]
Proposition~\ref{prop:binary-osculating} now gives that $f$ has a root
of multiplicity at least $r-k+1$.
\end{proof}

We obtain the binary case immediately.

\begin{corollary}\label{cor:binary-main}
If
\[
r>4p,
\]
then
\[
M_p(SL(V),S^rV)
=
\mathcal Q(T^pc_r)
=
T^pc_r.
\]
In particular, $M_p(SL(V),S^rV)$ is irreducible.
\end{corollary}

\begin{proof}
The inclusion
\[
T^pc_r
\subseteq
\mathcal Q(T^pc_r)
\]
holds by definition.

Conversely, let
\[
[f]\in\mathcal Q(T^pc_r).
\]
By Corollary~\ref{cor:binary-quadrics},
\[
\pi_m(f\cdot f)=0
\qquad
\text{for every }
p+1\leq m\leq\lfloor r/2\rfloor.
\]
Equivalently,
\[
(f,f)_{2m}=0
\qquad
\text{for every }
p+1\leq m\leq\lfloor r/2\rfloor.
\]
Applying Proposition~\ref{prop:weyman-binary} with $k=p+1$, and using
\[
r>4(p+1)-4=4p,
\]
we conclude that $f$ has a root of multiplicity at least
\[
r-(p+1)+1=r-p.
\]
Thus
\[
f=\ell^{r-p}h
\]
for some $0\neq\ell\in V$ and $0\neq h\in S^pV$. By Proposition~\ref{prop:binary-osculating},
\[
[f]\in T^pc_r.
\]
Therefore
\[
\mathcal Q(T^pc_r)=T^pc_r.
\]

Finally, $T^pc_r$ is irreducible, since by Proposition~\ref{prop:binary-osculating} it
is the image of the irreducible variety
\[
\mathbb P(V)\times\mathbb P(S^pV)
\]
under a projective morphism.
\end{proof}

\section{Proof of the main theorem}

We shall use the fact that the homogeneous ideal of a closed highest
weight orbit is generated by quadrics; see
\cite[Proposition~1.3.2]{GorodentsevKhoroshkinRudakov}.

We first record the following consequence.

\begin{lemma}\label{lem:intersection}
Let $d\geq p+2$. Then
\[
V_{d\lambda}\cap
v_\lambda^{\,d-p}S^p(V_\lambda)
=
E_p(d\lambda).
\]
\end{lemma}

\begin{proof}
Put $V=V_\lambda$, $v=v_\lambda$, and let
\[
I=I(X_\lambda)\subset S(V^*)
\]
be the homogeneous ideal of $X_\lambda$. Since
$V_{d\lambda}\subset S^dV$ is generated by the vectors $(gv)^d$,
$g\in G$, we have, with respect to the natural perfect pairing between
$S^d(V^*)$ and $S^dV$,
\begin{equation}\label{eq:cartan-annihilator}
I_d^\circ=V_{d\lambda}.
\end{equation}

Let $\mathfrak m_v\subset S(V^*)$ be the homogeneous ideal of the
point $[v]$. Under the embedding
\[
X_\lambda\longrightarrow\mathbb P(V_{d\lambda}),
\qquad
[g v]\longmapsto[(gv)^d],
\]
linear forms on $V_{d\lambda}$ are represented by degree-$d$ forms on
$V$, modulo $I_d$. We regard $E_p(d\lambda)^\circ$ here as the
annihilator of $E_p(d\lambda)\subset S^dV$ with respect to the natural
pairing between $S^d(V^*)$ and $S^dV$. The standard description of
osculating spaces in terms of jets shows that this annihilator consists
of the forms $F\in S^d(V^*)$ whose restriction to $X_\lambda$ vanishes
at $[v]$ to order at least $p+1$; see, for instance,
\cite[\S2]{MallavibarrenaPiene}. We claim that
\begin{equation}\label{eq:jet-annihilator}
E_p(d\lambda)^\circ
=
I_d+(\mathfrak m_v^{p+1})_d.
\end{equation}

One inclusion is immediate. For the converse, choose
$\alpha\in V^*$ with $\alpha(v)=1$ and dehomogenize on the affine
chart $\alpha\neq0$, taking $[v]$ as the origin. Let
$Q_1,\ldots,Q_s$ be quadratic generators of $I$. Since these
quadrics generate the ideal of $X_\lambda$ on this affine chart, if
$F$ vanishes on $X_\lambda$ to order at least $p+1$ at $[v]$, there
are polynomials $a_i$ such that
\[
\overline F-\sum_i a_i\overline{Q_i}
\]
vanishes at the origin to order at least $p+1$. Modulo terms of order
at least $p+1$, the $a_i$ may be chosen of degree at most $p$. Since
$d-2\geq p$, each $a_i$ can be homogenized to a form
$A_i\in S^{d-2}(V^*)$. It follows that
\[
F-\sum_i A_iQ_i
\in
(\mathfrak m_v^{p+1})_d,
\]
which proves \eqref{eq:jet-annihilator}.

Finally,
\begin{equation}\label{eq:point-annihilator}
\bigl((\mathfrak m_v^{p+1})_d\bigr)^\circ
=
v^{d-p}S^pV.
\end{equation}
Indeed, after choosing a decomposition $V=\mathbb Cv\oplus W$,
the degree-$d$ part of $\mathfrak m_v^{p+1}$ consists precisely of
the terms containing at least $p+1$ factors from $W^*$, and its
annihilator consists of the terms containing at most $p$ factors
from $W$.

Taking annihilators in \eqref{eq:jet-annihilator}, and using
\eqref{eq:cartan-annihilator} and \eqref{eq:point-annihilator}, gives
\[
E_p(d\lambda)
=
V_{d\lambda}\cap v^{d-p}S^pV.
\]
\end{proof}

\begin{theorem}\label{thm:main}
Let $G$ be a connected complex semisimple algebraic group, let $V_\lambda$
be an irreducible $G$-module of highest weight $\lambda$, and let $p\geq0$ and
$d\geq1$ be integers. If
\[
d>4p,
\]
then
\[
\mathcal Q(T^pX_{d\lambda})
=
T^pX_{d\lambda}.
\]
In particular, $\mathcal Q(T^pX_{d\lambda})$ is irreducible.
\end{theorem}

\begin{proof}
For $p=0$, the statement follows from the quadratic generation of
the ideal of $X_{d\lambda}$. We may therefore assume $p\geq1$.

Put
\[
V=V_\lambda,
\qquad
v=v_\lambda.
\]
If $\dim V=1$ there is nothing to prove, so assume $\dim V\geq2$.

The inclusion
\[
T^pX_{d\lambda}
\subseteq
\mathcal Q(T^pX_{d\lambda})
\]
holds by definition. Let
\[
[x]\in\mathcal Q(T^pX_{d\lambda}),
\qquad
0\neq x\in V_{d\lambda}.
\]

The action of $\mathfrak g$ on $S^dV$ satisfies the Leibniz rule.
Thus, for $\xi_1,\ldots,\xi_k\in\mathfrak g$ and $k\leq p$,
\[
\xi_1\cdots\xi_kv^d
\in
v^{d-k}S^kV
\subseteq
v^{d-p}S^pV.
\]
Consequently,
\[
E_p(d\lambda)\subseteq v^{d-p}S^pV,
\]
and, by $G$-equivariance,
\begin{equation}\label{eq:binary-containment}
gE_p(d\lambda)
\subseteq
(gv)^{d-p}S^pV
\qquad
\text{for every }g\in G.
\end{equation}

Let $\rho:V\to U$ be a surjective linear map with $\dim U=2$, and
let
\[
c_d(U)=\{[u^d]:0\neq u\in U\}
\subset\mathbb P(S^dU)
\]
be the rational normal curve. If
\[
q\in I(T^pc_d(U))_2,
\]
then the pull-back of $q$ by
\[
S^d\rho:S^dV\longrightarrow S^dU
\]
vanishes on $T^pX_{d\lambda}$. Indeed, by
\eqref{eq:binary-containment},
\[
S^d\rho\bigl(gE_p(d\lambda)\bigr)
\subseteq
\rho(gv)^{d-p}S^pU,
\]
which is either zero or the affine cone over a $p$-osculating space
of $c_d(U)$. Therefore the restriction to $V_{d\lambda}$ of the pull-back of
$q$ belongs to $I(T^pX_{d\lambda})_2$. Since
$[x]\in\mathcal Q(T^pX_{d\lambda})$, it follows that
\[
q(S^d\rho(x))=0.
\]
Hence, whenever $S^d\rho(x)\neq0$, the hypothesis $d>4p$ and
Corollary~\ref{cor:binary-main} give
\begin{equation}\label{eq:binary-section}
[S^d\rho(x)]\in T^pc_d(U).
\end{equation}

Regard $x\in S^dV$ as a homogeneous polynomial on $V^*$. If
$\dim V=2$, taking $\rho$ to be an isomorphism in
\eqref{eq:binary-section} gives
\[
x=\ell^{d-p}h
\]
for some $0\neq\ell\in V$ and $h\in S^pV$.

Assume $\dim V\geq3$ and factor
\[
x=P_1^{e_1}\cdots P_s^{e_s}
\]
into distinct irreducible homogeneous factors. Choose a general
projective line $L\subset\mathbb P(V^*)$ which is not contained in
$Z(x)$, avoids the singular loci and the pairwise intersections of
the hypersurfaces $P_i=0$, and meets each of them transversely. Such
lines form a non-empty open subset of the Grassmannian of projective
lines. Then every root of $x|_L$ lies on exactly one hypersurface
$P_i=0$ and has multiplicity equal to the corresponding exponent
$e_i$.

The line $L$ is of the form $\mathbb P(U^*)$ for a quotient
$\rho:V\to U$, and $x|_L=S^d\rho(x)\neq0$. By
\eqref{eq:binary-section}, $x|_L$ has a root of multiplicity at least
$d-p$. Hence
\[
e_i\geq d-p
\]
for some $i$. If $\deg P_i\geq2$, then
\[
d
\geq
e_i\deg P_i
\geq
2(d-p)
>
d,
\]
a contradiction. Therefore $P_i$ is linear. Thus, in all cases,
\begin{equation}\label{eq:large-linear-factor}
x=\ell^{d-p}h.
\end{equation}
for some $0\neq\ell\in V$ and $h\in S^pV$.

Write more precisely
\[
x=\ell^a h_0,
\qquad
\ell\nmid h_0,
\]
where $a$ is the multiplicity of the linear factor $\ell$ in $x$.
By \eqref{eq:large-linear-factor},
\[
a\geq d-p.
\]
Since $p\geq1$ and $d>4p$, in particular $a\geq2$. We claim that
$[\ell]\in X_\lambda$.

Suppose otherwise. Since the ideal of $X_\lambda$ is generated by
quadrics, there exists
\[
q\in I(X_\lambda)_2
\]
such that
\[
q(\ell)\neq0.
\]

Since $\ell\nmid h_0$, the image $\overline h_0$ of $h_0$ in
\[
S^{d-a}(V/\mathbb C\ell)
\]
is non-zero. Choose
\[
\overline\varphi\in
S^{d-a}((V/\mathbb C\ell)^*)
\]
with $\overline\varphi(\overline h_0)\neq0$, and let
$\varphi\in S^{d-a}(V^*)$ be its pull-back. Thus
$\varphi(h_0)=\overline\varphi(\overline h_0)\neq0$. When $a=d$,
this simply means that $\varphi$ is a non-zero constant. Choose also
$\alpha\in V^*$ with $\alpha(\ell)=1$.

Since
\[
q\,\alpha^{a-2}\varphi\in I(X_\lambda)_d
\]
and $x\in V_{d\lambda}$, equation~\eqref{eq:cartan-annihilator}
gives
\begin{equation}\label{eq:contradiction-zero}
(q\,\alpha^{a-2}\varphi)(x)=0.
\end{equation}

On the other hand, if $a<d$, viewed as a symmetric multilinear
form, $\varphi$ vanishes whenever one of its arguments lies in
$\mathbb C\ell$. Thus, when $q\,\alpha^{a-2}\varphi$ is evaluated
on $x=\ell^a h_0$, the $d-a$ arguments of $\varphi$ must all come
from $h_0$. The same conclusion is vacuous when $a=d$. Therefore
\[
(q\,\alpha^{a-2}\varphi)(x)
=
c\,q(\ell)\,\varphi(h_0)
\]
for some non-zero constant $c$. Since
$q(\ell)\neq0$ and $\varphi(h_0)\neq0$, this contradicts
\eqref{eq:contradiction-zero}. Hence
\[
[\ell]\in X_\lambda.
\]

Since $T^pX_{d\lambda}$ and its quadratic hull are $G$-stable, after
acting by an element of $G$ and absorbing a non-zero scalar in $h_0$,
we may assume $\ell=v$. Since $a\geq d-p$,
\[
x\in
V_{d\lambda}\cap v^{d-p}S^pV.
\]
The hypothesis $d>4p$ implies $d\geq p+2$, so
Lemma~\ref{lem:intersection} gives
\[
x\in E_p(d\lambda).
\]
Therefore
\[
[x]\in T^p_{[v^d]}X_{d\lambda}
\subseteq
T^pX_{d\lambda}.
\]
Undoing the action of $G$, we conclude that
\[
\mathcal Q(T^pX_{d\lambda})
\subseteq
T^pX_{d\lambda}.
\]
Hence equality holds.

Finally, by Proposition~\ref{prop:pbw-osculating},
\[
T^pX_{d\lambda}
=
\overline{G\cdot\mathbb P(E_p(d\lambda))}.
\]
Since $G$ and $\mathbb P(E_p(d\lambda))$ are irreducible, the
closure of their image is irreducible. Thus $T^pX_{d\lambda}$, and
therefore its quadratic hull, is irreducible.
\end{proof}

\section{Examples and limits of the stable range}

The main theorem gives
\[
M_p(G,V_{d\lambda})
=
T^pX_{d\lambda}
\]
whenever $d>4p$. We begin with two examples outside this range.
The first shows that the quadratic hull may remain irreducible while
being strictly larger than the osculating variety. The second shows
that, outside the stable range, the quadratic hull may acquire
additional irreducible components. We then finish with a related
example, outside the osculating setting, in which a reducible cubic
base locus gives rise to a reducible quadratic hull.

\subsection{An irreducible binary quadratic hull outside the stable range}

The stable range in the binary case guarantees that the quadratic
hull of the osculating variety coincides with the osculating variety
itself. Outside this range, the quadratic hull may still be
irreducible without coinciding with the osculating variety. The
following example illustrates this phenomenon.

\begin{proposition}\label{prop:binary-example}
Let $V=\mathbb C^2$. Then
\[
M_2(SL(V),S^7V)=\mathcal Q(T^2c_7)
\]
is irreducible of dimension $4$. In particular,
\[
T^2c_7\subsetneq M_2(SL(V),S^7V).
\]
\end{proposition}

\begin{proof}
Choose coordinates $u,v$ on $V$ and write
\[
f=\sum_{i=0}^7 a_i u^{7-i}v^i.
\]
By Corollary~\ref{cor:binary-quadrics}, only the summand corresponding
to $m=3$ occurs in $I(T^2c_7)_2$. Using the explicit formulas of
\cite[Corollary~2.2]{Massri2013}, the variety
\[
M_2(SL(V),S^7V)
\]
is therefore defined by the three quadrics
\[
\begin{aligned}
q_0&=
a_0a_6-6a_1a_5+15a_2a_4-10a_3^2,\\
q_1&=
a_0a_7-5a_1a_6+9a_2a_5-5a_3a_4,\\
q_2&=
a_1a_7-6a_2a_6+15a_3a_5-10a_4^2.
\end{aligned}
\]

We work in the affine chart $a_0=1$. The equations $q_0=q_1=0$
determine $a_6$ and $a_7$:
\[
a_6
=
6a_1a_5-15a_2a_4+10a_3^2
\]
and
\[
a_7
=
5a_1a_6-9a_2a_5+5a_3a_4.
\]
After substituting these expressions in $q_2$, the remaining
equation has the form
\[
Aa_5+B=0,
\]
where
\[
A
=
3(2a_1^3-3a_1a_2+a_3)
\]
and
\[
B=
-15a_1^2a_2a_4
+10a_1^2a_3^2
+a_1a_3a_4+18a_2^2a_4
-12a_2a_3^2
-2a_4^2.
\]

We claim that $Aa_5+B$ is irreducible. The polynomial $A$ is
irreducible in $\mathbb C[a_1,a_2,a_3,a_4]$, since, up to a non-zero
constant, it is monic of degree one in $a_3$. Moreover,
\[
\gcd(A,B)=1.
\]
Indeed, regard $B$ as a polynomial in $a_4$. Since $A$ is independent
of $a_4$, if $A$ divided $B$, it would divide every coefficient of
$B$ as a polynomial in $a_4$. This is impossible because the
coefficient of $a_4^2$ is the non-zero constant $-2$.

Thus $Aa_5+B$ is primitive as a polynomial in $a_5$ over
$\mathbb C[a_1,a_2,a_3,a_4]$. Over the corresponding fraction field
it has degree one, and hence it is irreducible. By Gauss' lemma it is
irreducible in $\mathbb C[a_1,\ldots,a_5]$. Consequently,
\[
M_2(SL(V),S^7V)\cap\{a_0\neq0\}
\]
is isomorphic to the irreducible hypersurface
\[
Z(Aa_5+B)\subset\mathbb A^5_{a_1,\ldots,a_5}.
\]
In particular, this open subset has dimension $4$.

We now prove that the whole projective variety is irreducible. Since
$SL(V)$ is connected, each irreducible component of
$M_2(SL(V),S^7V)$ is $SL(V)$-invariant. Fix the Borel subgroup for
which $u^7$ is a highest weight vector. By the Borel fixed point
theorem, every irreducible component contains a Borel fixed point.
The only Borel fixed point in $\mathbb P(S^7V)$ is
\[
[u^7]=[1:0:\cdots:0],
\]
so every irreducible component meets the chart $a_0\neq0$.

Let $X_1,\ldots,X_s$ be the irreducible components and put
\[
U=M_2(SL(V),S^7V)\cap\{a_0\neq0\}.
\]
We have just proved that $U$ is irreducible and that
$X_i\cap U\neq\varnothing$ for every $i$. Since
\[
U=\bigcup_{i=1}^s(X_i\cap U)
\]
and each $X_i\cap U$ is closed in $U$, irreducibility of $U$ implies
that $U\subseteq X_j$ for some $j$. For every $i$, the set
$X_i\cap U$ is a non-empty open subset of $X_i$, hence it is dense in
$X_i$. Since $X_i\cap U\subseteq X_j$, taking closures gives
$X_i\subseteq X_j$. Thus all the irreducible components coincide,
and $M_2(SL(V),S^7V)$ is irreducible.

Finally, Proposition~\ref{prop:binary-osculating} gives a
surjective morphism
\[
\mathbb P(V)\times\mathbb P(S^2V)
\longrightarrow
T^2c_7.
\]
For a general form $\ell^5h$, with $\ell\nmid h$, the factor $\ell$
is uniquely determined. Hence this morphism is generically one-to-one and
\[
\dim T^2c_7=1+2=3.
\]
On the other hand,
\[
\dim M_2(SL(V),S^7V)=4.
\]
Since
\[
T^2c_7\subset M_2(SL(V),S^7V),
\]
the inclusion is strict.
\end{proof}

\begin{remark}
This example shows that the equality
\[
\mathcal Q(T^pc_r)=T^pc_r
\]
is stronger than the irreducibility of the quadratic hull. In the
example above the quadratic hull is irreducible, but it is strictly
larger than the osculating variety.

The calculation also suggests that the stable bound $r>4p$ is not
expected to be optimal for irreducibility in the binary case. The
same elimination procedure can be carried out in further examples,
and it is natural to ask for the full range of pairs $(r,p)$ for
which $M_p(SL(V),S^rV)$ is irreducible, even when it does not coincide
with $T^pc_r$. A systematic analysis of this question seems to
require a separate study of the corresponding systems of quadratic
equations, and we leave it for future work.
\end{remark}

\subsection{A reducible quadratic hull for a restricted chordal variety}

Let $V=\mathbb C^N$ and consider the Pl\"ucker embedding
\[
\operatorname{Gr}_k(V)\subset\mathbb P(\wedge^kV).
\]
Following the notation of \cite[Hard Exercise~15.44]{FultonHarris},
write the multiplicity-free decomposition
\[
S^2(\wedge^kV)=\bigoplus_{i\geq0}\Theta_{2i},
\]
where $\Theta_0$ is the Cartan component of highest weight
$2\omega_k$. For $q\geq1$, let
\[
C^q(\operatorname{Gr}_k(V))
\]
denote the $q$-restricted chordal variety, namely the closure of the
union of the lines joining two points corresponding to $k$-planes
$L,M\subset V$ such that
\[
\dim(L\cap M)\geq k-2q+1.
\]
As stated in \cite[Hard Exercise~15.44]{FultonHarris},
\begin{equation}\label{eq:restricted-chordal-quadrics}
I(C^q(\operatorname{Gr}_k(V)))_2
=
\bigoplus_{i\geq q}\Theta_{2i}^*.
\end{equation}

The following example identifies one of these quadratic hulls with
one of the varieties considered in this article.

\begin{proposition}\label{prop:grassmannian-reducible}
Let $V=\mathbb C^{12}$. Then
\[
M_1(SL(V),\wedge^4V)
=
\mathcal Q(C^2(\operatorname{Gr}_4(V)))
=
\left\{
[\omega]\in\mathbb P(\wedge^4V):
\omega\wedge\omega=0
\right\}.
\]
Moreover, this variety is reducible.
\end{proposition}

\begin{proof}
Here the closed highest weight orbit in
$\mathbb P(\wedge^4V)$ is the Pl\"ucker Grassmannian
$\operatorname{Gr}_4(V)$. Let
\[
v=e_1\wedge e_2\wedge e_3\wedge e_4
\]
be a highest weight vector and put
\[
A=\langle e_1,e_2,e_3,e_4\rangle.
\]
By the definition of the PBW filtration,
\[
E_1(\omega_4)=\wedge^3A\wedge V.
\]

For $k=4$ the symmetric square decomposes as
\[
S^2(\wedge^4V)
=
\Theta_0\oplus\Theta_2\oplus\Theta_4,
\]
and
\[
\Theta_4\simeq\wedge^8V.
\]
Up to a non-zero scalar, the projection onto $\Theta_4$ is the
exterior product
\[
S^2(\wedge^4V)\longrightarrow\wedge^8V,
\qquad
x\cdot y\longmapsto x\wedge y.
\]

If $x,y\in E_1(\omega_4)$, each decomposable summand of $x$ and $y$
contains three factors from the $4$-dimensional space $A$.
Consequently, the exterior product of any two such summands contains
six factors from $A$ and therefore vanishes. Hence
\[
x\wedge y=0
\]
for every $x,y\in E_1(\omega_4)$, and it follows that
\[
F_1(\omega_4)
\subseteq
\Theta_0\oplus\Theta_2.
\]

Since $v\cdot v\in S^2(E_1(\omega_4))$, its $SL(V)$-span is the
Cartan component, and therefore
\[
\Theta_0\subseteq F_1(\omega_4).
\]
We claim that $\Theta_2$ also occurs. Consider
\[
x=e_{1234}+e_{1235}+e_{1246}\in E_1(\omega_4),
\]
where $e_{i_1\cdots i_r}$ denotes
$e_{i_1}\wedge\cdots\wedge e_{i_r}$. Since all three summands
contain $e_1\wedge e_2$, we have
\[
x\wedge x=0,
\]
so the $\Theta_4$-component of $x\cdot x$ vanishes. On the other
hand, $x$ is not decomposable. Indeed, the Pl\"ucker relation
\[
p_{1234}p_{1256}
-
p_{1235}p_{1246}
+
p_{1236}p_{1245}
=0
\]
does not vanish at $x$.

If the $\Theta_2$-component of $x\cdot x$ were also zero, then
$x\cdot x$ would belong to $\Theta_0=F_0(\omega_4)$. By
Proposition~\ref{prop:quadrics}, this would imply that all the
quadratic equations of $\operatorname{Gr}_4(V)$ vanish at $x$, in
particular the Pl\"ucker relation displayed above, a contradiction.
Thus $\Theta_2$ occurs in $F_1(\omega_4)$. Since $F_1(\omega_4)$ is
$SL(V)$-stable and the above decomposition is multiplicity-free,
\[
F_1(\omega_4)=\Theta_0\oplus\Theta_2.
\]

By Proposition~\ref{prop:quadrics},
\[
I(T^1\operatorname{Gr}_4(V))_2
=
\Theta_4^*.
\]
On the other hand, taking $q=2$ in
\eqref{eq:restricted-chordal-quadrics} gives
\[
I(C^2(\operatorname{Gr}_4(V)))_2
=
\Theta_4^*.
\]
Consequently,
\[
M_1(SL(V),\wedge^4V)
=
\mathcal Q(C^2(\operatorname{Gr}_4(V))).
\]
Since the projection onto $\Theta_4\simeq\wedge^8V$ is the exterior
product, this common quadratic hull is
\[
M_1(SL(V),\wedge^4V)
=
\left\{
[\omega]\in\mathbb P(\wedge^4V):
\omega\wedge\omega=0
\right\}.
\]

It remains to prove that this variety is reducible. Consider the
projective bundle over $\mathbb P(V)$ whose fiber over $[\ell]$ is
\[
\mathbb P\left(\wedge^3(V/\mathbb C\ell)\right).
\]
The natural map to $\mathbb P(\wedge^4V)$ sends a point over $[\ell]$
to the class of a $4$-form divisible by $\ell$. Let $Z$ denote its
image. Thus
\[
Z=
\left\{
[\ell\wedge\eta]:
0\neq\ell\in V,\quad \ell\wedge\eta\neq0
\right\}.
\]
Since the source is projective and irreducible, $Z$ is closed and
irreducible. Moreover,
\[
Z\subset M_1(SL(V),\wedge^4V).
\]

Let
\[
\omega=\ell\wedge\eta
\]
be general in $Z$ and choose a decomposition
\[
V=\mathbb C\ell\oplus U,
\qquad
\dim U=11,
\]
with $\eta\in\wedge^3U$. By
\cite[Corollary~3.2 and Proposition~4.5]{LundqvistNicklasson}, for a
general $3$-form $\eta$ in $11$ variables the map
\[
\eta\wedge-:
U\longrightarrow\wedge^4U
\]
is injective, and
\[
\dim\ker
\left(
\eta\wedge-:
\wedge^4U\longrightarrow\wedge^7U
\right)
=11.
\]
The injectivity of the first map implies
\[
\{w\in V:w\wedge\omega=0\}
=
\mathbb C\ell.
\]
Thus the linear factor $[\ell]$ is uniquely determined by a general
point of $Z$. The projective bundle above has dimension
\[
11+\left(\binom{11}{3}-1\right)=175,
\]
and the uniqueness of $[\ell]$ shows that its map onto $Z$ is generically
one-to-one. Hence
\begin{equation}\label{eq:dim-Z}
\dim Z=175.
\end{equation}

We now compute the tangent space to
$M_1(SL(V),\wedge^4V)$ at a general point of $Z$. Consider
\[
\Phi:\wedge^4V\longrightarrow\wedge^8V,
\qquad
\xi\longmapsto\xi\wedge\xi.
\]
Its differential at $\omega=\ell\wedge\eta$ is
\[
d\Phi_\omega(\xi)=2\omega\wedge\xi.
\]
Writing
\[
\xi=\ell\wedge\alpha+\beta,
\qquad
\alpha\in\wedge^3U,
\qquad
\beta\in\wedge^4U,
\]
we obtain
\[
d\Phi_\omega(\xi)
=
2\ell\wedge\eta\wedge\beta.
\]
The differentials of the defining quadrics therefore cut out
$\ker(d\Phi_\omega)$. Hence the Zariski tangent space of the reduced
zero locus satisfies
\[
\dim T_{[\omega]}M_1(SL(V),\wedge^4V)
\leq
\binom{11}{3}
+
\dim\ker
\left(
\eta\wedge-:
\wedge^4U\longrightarrow\wedge^7U
\right)
-1
=
175.
\]

Let $C$ be an irreducible component of
$M_1(SL(V),\wedge^4V)$ containing $Z$. Then
\[
175
=
\dim Z
\leq
\dim C
\leq
\dim T_{[\omega]}M_1(SL(V),\wedge^4V)
\leq
175.
\]
Thus $\dim C=\dim Z$, and since $Z\subseteq C$ is closed and
irreducible, we obtain $C=Z$. Hence $Z$ is an irreducible component
of $M_1(SL(V),\wedge^4V)$.

It remains to show that $Z$ is not the whole variety. Let
\[
U_0=\langle e_1,\ldots,e_7\rangle\subset V
\]
and consider
\[
\psi
=
e_{3456}+e_{1256}+e_{1234}
+e_{2467}-e_{2357}-e_{1457}-e_{1367}
\in\wedge^4U_0.
\]
Since $\dim U_0=7$,
\[
\psi\wedge\psi=0,
\]
and hence
\[
[\psi]\in M_1(SL(V),\wedge^4V).
\]

We claim that $\psi$ is not divisible by a non-zero vector. A direct
computation gives
\[
\begin{aligned}
e_1\wedge\psi
 &=
 e_{13456}+e_{12467}-e_{12357},\\
e_2\wedge\psi
 &=
 e_{23456}+e_{12457}+e_{12367},\\
e_3\wedge\psi
 &=
 e_{12356}-e_{23467}+e_{13457},\\
e_4\wedge\psi
 &=
 e_{12456}-e_{23457}-e_{13467},\\
e_5\wedge\psi
 &=
 e_{12345}+e_{24567}-e_{13567},\\
e_6\wedge\psi
 &=
 e_{12346}+e_{23567}+e_{14567},\\
e_7\wedge\psi
 &=
 e_{34567}+e_{12567}+e_{12347}.
\end{aligned}
\]
The supports of these seven vectors in the standard basis of
$\wedge^5U_0$ are pairwise disjoint. Hence they are linearly
independent, and
\[
w\wedge\psi=0,
\qquad
w\in U_0,
\]
implies $w=0$.

Choose a complement $V=U_0\oplus U_1$. If
$w=w_0+w_1$, with $w_i\in U_i$, satisfies
\[
w\wedge\psi=0,
\]
then
\[
w_0\wedge\psi\in\wedge^5U_0,
\qquad
w_1\wedge\psi\in U_1\wedge\wedge^4U_0.
\]
These belong to distinct summands of $\wedge^5V$, so they vanish
separately. The preceding computation gives $w_0=0$, while
$w_1\wedge\psi=0$ implies $w_1=0$, since under the natural identification
$U_1\wedge\wedge^4U_0\simeq U_1\otimes\wedge^4U_0$ this map is
$w_1\mapsto w_1\otimes\psi$, and $\psi\neq0$. Thus
\[
w\wedge\psi=0
\quad\Longrightarrow\quad
w=0.
\]
In particular, $\psi$ cannot be written as
\[
\psi=w\wedge\gamma
\]
with $w\neq0$, and therefore
\[
[\psi]\notin Z.
\]

Since $Z$ is an irreducible component and
\[
Z\neq M_1(SL(V),\wedge^4V),
\]
the variety $M_1(SL(V),\wedge^4V)$ is reducible.
\end{proof}

\begin{remark}
The Pl\"ucker representation has highest weight $\omega_4$, so in
the notation of the main theorem it corresponds to $d=1$. Thus the
preceding example lies outside the stable range $d>4p$. It gives a
concrete instance in which the quadratic hull associated with the
osculating construction acquires additional irreducible components.
In particular, in this case the actual zero locus of the quadrics
appearing in \cite[Hard Exercise~15.44]{FultonHarris} is reducible.
\end{remark}

\subsection{A reducible quadratic hull from the trifocal variety}

The preceding two examples concern osculating varieties outside the
stable range. We finish with a related phenomenon which does not come
from the osculating construction. The following observation turns a
reducible cubic base locus into a reducible quadratic hull after a
second Veronese embedding.

\begin{lemma}\label{lem:veronese-quadratic-hull}
Let
\[
Y\subset\mathbb P(W)
\]
be a projective variety and let
\[
\nu_2:\mathbb P(W)\longrightarrow\mathbb P(S^2W)
\]
be the second Veronese embedding. Assume that
\begin{equation}\label{eq:veronese-degree-four}
I(Y)_4=W^*I(Y)_3.
\end{equation}
Then
\[
\mathcal Q(\nu_2(Y))
=
\nu_2\bigl(Z(I(Y)_3)\bigr).
\]
In particular, if $Z(I(Y)_3)$ is reducible, then
$\mathcal Q(\nu_2(Y))$ is reducible.
\end{lemma}

\begin{proof}
The multiplication map
\[
S^2(S^2W^*)
\longrightarrow
S^4W^*
\]
is surjective, and its kernel is
\[
I(\nu_2(\mathbb P(W)))_2.
\]
Since
\[
I(\nu_2(\mathbb P(W)))_2
\subseteq
I(\nu_2(Y))_2,
\]
all the quadratic equations of the Veronese variety occur among the
quadratic equations of $\nu_2(Y)$. Therefore
\[
\mathcal Q(\nu_2(Y))
\subseteq
\nu_2(\mathbb P(W)).
\]
Inside the Veronese variety, the remaining quadratic equations pull
back exactly to $I(Y)_4$. Hence
\[
\mathcal Q(\nu_2(Y))
=
\nu_2\bigl(Z(I(Y)_4)\bigr).
\]

By \eqref{eq:veronese-degree-four},
\[
Z(I(Y)_4)
=
Z(W^*I(Y)_3)
=
Z(I(Y)_3).
\]
Indeed, one inclusion is immediate. Conversely, if
$x\notin Z(I(Y)_3)$, choose $f\in I(Y)_3$ with $f(x)\neq0$ and
$\ell\in W^*$ with $\ell(x)\neq0$. Then
\[
\ell f\in I(Y)_4
\]
and
\[
(\ell f)(x)\neq0.
\]
This proves the equality. Since $\nu_2$ is a closed embedding, it
preserves irreducible components, and the final assertion follows.
\end{proof}

We apply the lemma to the trifocal variety
\[
X\subset
\mathbb P(A^*\otimes B^*\otimes C),
\qquad
\dim A=\dim B=\dim C=3,
\]
and put
\[
W=A^*\otimes B^*\otimes C.
\]
Aholt and Oeding prove in
\cite[Theorem~1.1]{AholtOeding} that the prime ideal $I(X)$ is
minimally generated by $10$ polynomials of degree $3$, $81$
polynomials of degree $5$, and $1980$ polynomials of degree $6$.
In particular, there are no minimal generators in degree $4$, and
hence
\[
I(X)_4=W^*I(X)_3.
\]
Moreover, \cite[Proposition~7.3]{AholtOeding} identifies the zero
locus of these ten cubic equations as a union of four irreducible
components, one of which is the trifocal variety itself.

Lemma~\ref{lem:veronese-quadratic-hull} therefore gives
\[
\mathcal Q(\nu_2(X))
=
\nu_2\bigl(Z(I(X)_3)\bigr).
\]
Consequently,
\[
\mathcal Q(\nu_2(X))
\]
has four irreducible components and, in particular, is reducible.

The first two examples of this section show that outside the stable
range the quadratic hull of an osculating variety may either remain
irreducible while becoming strictly larger than the osculating
variety, or acquire additional irreducible components. The trifocal
example shows that reducible quadratic hulls also arise naturally
beyond the osculating setting.

\bibliographystyle{plain}
\bibliography{citas2}

\end{document}